\documentclass[10pt]{article}
\usepackage{a4}
\usepackage{latexsym}
\usepackage{amsfonts}
\usepackage{amssymb}
\usepackage{amsmath}
\usepackage{amsthm}
\newtheorem{theorem}{Theorem}[section]
\newtheorem{lemma}[theorem]{Lemma}

\newtheorem{application}[theorem]{Application}
\title{An improved sum-product estimate for general finite fields}
\author{Oliver Roche-Newton}
\date{}

\begin{document}
\maketitle
\begin{abstract}

This paper improves on a sum-product estimate obtained by Katz and Shen for subsets of a finite field whose order is not prime.

\end{abstract}

\section{Introduction}

Let $A$ be a subset of $F=\mathbb{F}_{p^n}$ for some prime $p$ and some $n\in{\mathbb{N}}$. Consider the sumset and productset of $A$, defined respectively as

$$A+A:=\{a+b:a,b\in{A}\},$$
$$A\cdot{A}:=\{ab:a,b\in{A}\}.$$

An interesting problem is to establish lower bounds on the quantity $\max{\{|A+A|,|A\cdot{A}|\}}$. The existence of non-trivial bounds was first establish by Bourgain, Katz and Tao\cite{BKT}. Garaev\cite{garaev} established the first quantitative sum-product estimate for fields of prime order. Garaev's result can be stated as follows:

\begin{theorem}\label{theorem:garaev} Let $A\subset{\mathbb{F}_p}$ for some prime $p$, such that $|A|\leq{p^{7/13}}(\log{p})^{-4/13}$. Then

$$\max\{|A+A|,|A\cdot{A}|\}\gg{\frac{|A|^{15/14}}{(\log{|A|})^{2/7}}}.$$

\end{theorem}

This result was generalised, for finite fields whose order is not necessarily prime, by Katz and Shen\cite{KS} in the form of the following theorem.

\begin{theorem}\label{theorem:KStheorem}Let $F=\mathbb{F}_{p^n}$ be a finite field. Suppose that $A$ is a subset of $F$ so that  for any subfield $G\subseteq{F}$, and any elements $c,d\in{F}$,

\begin{equation}
|A\cap{(cG+d)}|\leq{\max\{|G|^{1/2},|A|^{18/19}\}}.
\label{condition1}
\end{equation}

Then it must be the case that

$$\max{\{|A+A|,|A\cdot{A}|\}}\gg{\frac{|A|^{20/19}}{(\log{|A|})^{\alpha}}},$$

where $0<{\alpha}\leq{1}$ is some absolute constant.

\end{theorem}

Since Garaev's sum-product estimate\cite{garaev}, there have been a number of small improvements made courtesy of more subtle arguments by Katz-Shen\cite{KSprime}, Bourgain-Garaev\cite{BG}, Shen\cite{ShenSP}, Li\cite{Li} and most recently Rudnev\cite{mishaSP}; so the following result of Rudnev represents the state of the art.

\begin{theorem}Let $A\subset{F_p^*}$ with $|A|<p^{1/2}$ and $p$ large. Then

$$\max{\{|A+A|,|A\cdot{A}|\}}\gg{\frac{|A|^{12/11}}{(\log{|A|})^{4/11}}}.$$

\end{theorem}

This paper seeks to update the result of Katz and Shen for general finite fields, by using the refinements that allowed for the improved sum-product estimates in the prime setting in \cite{KSprime}, \cite{BG}, \cite{ShenSP} and \cite{mishaSP}. As well as adopting the techniques from the more recent sum-product estimates in prime fields, a slightly different argument is made in what will later be called case 3, which leads to a further improvement in the final estimate. The outcome is the following result:

\begin{theorem}\label{theorem:main} 

Let $F=F_{p^n}$ be a finite field and suppose that $A$ is a subset of $F^*$ with the following property. For any subfield $G\subset{F}$ and any elements $c,d\in{F}$,

\begin{equation}
|A\cap{(cG+d)}|\leq{\max\{|G|^{14/25},|A|^{6/7}\}}.
\label{condition2}
\end{equation}

Then we get the sum-product inequality

$$\max\{|A+A|,|A\cdot{A}|\}\gg{\frac{|A|^{\frac{15}{14}}}{\log|A|}}.$$

\end{theorem}

Note that the final estimate here is aligned with the original quantitative bound Garaev obtained in Theorem \ref{theorem:garaev}.

\subsection{Some remarks concerning Theorem \ref{theorem:main}}

As well as improving the exponent in the final bound, the first term from condition \eqref{condition1} has been increased in \eqref{condition2}, thus relaxing one of the conditions. However, it seems necessary to also tighten one of the conditions, as the second term from \eqref{condition1} becomes smaller in \eqref{condition2}.

The statement of Theorem \ref{theorem:main} is in fact quite flexible. The tools used do not distinguish between addition and subtraction, which means that the difference set, $A-A$, can replace $A+A$ in the above. It is also possible to get a sum-ratio estimate, where $A\cdot{A}$ is replaced by the ratio set $A:A$, which is defined as follows:

$$A:A\stackrel{\mathrm{def}}{=}\left\{\frac{a}{b}\text{ such that }a,b\in{A},b\neq{0}\right\}.$$

The presentation of a sum-ratio estimate proof is of a very similar nature, but rather more simple, as Rudnev \cite{mishaSP} alluded to. Lemma \ref{theorem:misha} can be replaced by a more straightforward pigeonholing argument and there is no need for dyadic pigeonholing, which means less technicalities are required and no logarithmic factor appears in the final estimate. One more benefit of this simplified argument is that the condition from Theorem \ref{theorem:KStheorem} which was made stricter for the main result in this paper (the second term in \eqref{condition1} and \eqref{condition2}), can be relaxed slightly.

\subsection{Notation}

Throughout this paper, the symbols $\ll,\gg$ and $\approx$ are used to suppress constants. For example, $X\ll{Y}$ means that there exists some absolute constant $C$ such that $X<CY$. $X\approx{Y}$ means that $X\ll{Y}$ and $Y\ll{X}$. When describing such a rough inequality in general language, inverted commas are used in an effort to avoid confusion. For example, the statement that ``$X$ is `at most' $Y$'' tells us that $X\ll{Y}$. A similar meaning is attached to `at least'.

\subsection{Acknowledgements} 

The author would like to thank Timothy Jones and Misha Rudnev for some extremely helpful conversations and typographical corrections.

\section{Preliminary results} 

Before proving the main theorem, it is necessary to state some previous results. The first is a lemma from the paper of Katz and Shen\cite{KS}.

\begin{lemma}\label{theorem:KatzShen} Let $A\subset{F}$ and let $G\subseteq{F}$ be a subfield such that

$$\frac{A-A}{A-A}\subseteq{G}.$$

Then there exist $c,d\in{F}$ such that

$$A\subseteq{cG+d}.$$

\end{lemma}

The next two results are well-known within arithmetic combinatorics, and have been crucial to all known quantitative sum-product estimates over finite fields. The first is due to Pl\"{u}nnecke and Rusza (see \cite{ruzsa}), whilst the second is a generalisation which Katz and Shen\cite{KSprime} used to obtain an improve the original quantitative sum-product estimate of Garaev\cite{garaev}.

\begin{lemma}\label{theorem:SPlun}Let $X,B_1,...,B_k$ be subsets of $F$. Then

$$|B_1+...+B_k|\leq{\frac{|X+B_1|...|X+B_k|}{|X|^{k-1}}}.$$

\end{lemma}

\begin{lemma}\label{theorem:CPlun}Let $X,B_1,...,B_k$ be subsets of $F$. Then for any $0<\epsilon<1$, there exists a subset $X'\subseteq{X}$, with $|X'|\geq{(1-\epsilon)|X|}$, and some constant $C(\epsilon)$, such that

$$|X'+B_1+...+B_k|\leq{C(\epsilon)\frac{|X+B_1|...|X+B_k|}{|X|^{k-1}}}.$$

\end{lemma}

Finally, we will need the following covering lemma, which appeared in sum-product estimates for the first time in \cite{ShenSP}. An application of this lemma has been the key to the two of the most recent improvements to the sum-product estimate over prime fields (see \cite{mishaSP}, \cite{ShenSP},).

\begin{lemma}\label{theorem:covering} Let $X$ and $Y$ be additive sets. Then for any $\epsilon\in{(0,1)}$ there is some constant $C(\epsilon)$, such that at least $(1-\epsilon)|X|$ of the elements of $X$ can be covered by $C(\epsilon)\frac{\min\{|X+Y|,|X-Y|\}}{|Y|}$ translates of $Y$.

\end{lemma}

The rest of this paper is devoted to proving Theorem \ref{theorem:main}.

\section{Proof of Theorem \ref{theorem:main}}

Let $A$ be a set satisfying the conditions of Theorem \ref{theorem:main}, and suppose that $|A+A|,|A\cdot{A}|\leq{K|A|}$. The aim is to show that $K\gg{\frac{|A|^{1/14}}{\log{|A|}}}$.

At the outset, apply Lemma \ref{theorem:CPlun} to identify some subset $A'\subset{A}$, with cardinality $|A'|\approx{|A|}$, so that

$$|A'+A'+A'+A'|\ll{\frac{|A+A|^3}{|A|^2}}\ll{K^3|A|}.$$

Next apply Lemma \ref{theorem:CPlun} again to identify a further subset $A''$, with cardinality $|A''|\approx{|A|}$, such that

$$|A''+A''+A''|\ll{\frac{|A+A|^2}{|A|}}\ll{K^2|A|}.$$

Since many more refinements of $A$ are needed throughout the proof, this first change is made without a change in notation. So, throughout the rest of the proof, when the set $A$ is referred to, we are really talking about the large subset $A''$, which satisfies each of the above inequalities. In other words, we assume that

\begin{equation}
|A+A+A+A|\ll{K^3|A|}.
\label{finaltouch}
\end{equation}

and

\begin{equation}
|A+A+A|\ll{K^2|A|}
\label{finaltouch2}
\end{equation}

hold, for our set $A$.

Consider the point set $A\times{A}\subset{F\times{F}}$. The multiplicative energy, $E(A)$, of $A$ is defined to be the number of solutions to

\begin{equation}
\frac{a_1}{a_2}=\frac{a_3}{a_4},
\label{menergy}
\end{equation}

such that $a_1,a_2,a_3,a_4\in{A}$. Let $\mathcal{L}$ be the set of all lines through the origin. By definition,

$$E(A)=\sum_{l\in{\mathcal{L}}}|l\cap({A\times{A}})|^2.$$

By dyadic pigeonholing, a technique familiar from most known sum-product estimates (see \cite{garaev},\cite{KSprime},\cite{KS},\cite{ShenSP},\cite{mishaSP}), we can identify a popular dyadic group. This is done by partitioning lines through the origin according to their popularity. Speaking more precisely, let $\mathcal{L}_i$ be the set of all lines $l$ through the origin such that $2^i\leq{|l\cap{A\times{A}}|}<{|2^{i+1}|}$. Therefore the multiplicative energy may be rewritten in the form

$$E(A)=\sum_{j=0}^{\log{|A|}}\sum_{l\in{\mathcal{L}_j}}|l\cap({A\times{A}})|.$$

Then, by elementary piegeonholing, one may choose a particular dyadic group, which contributes more than the average for this sum. Therefore, there exists a set of $L$ lines, each supporting $\approx{N}$ points from $A\times{A}$ such that

\begin{equation}
M:=LN^2\gg{\frac{E(A)}{\log{|A|}}}\geq{\frac{|A|^4}{|A\cdot{A}|\log{|A|}}}\geq{\frac{|A|^3}{K\log{|A|}}},
\label{dyadic}
\end{equation}

where the second inequality comes from the standard Cauchy-Schwarz bound on the multiplicative energy (see \cite{TV}). Refine the point set by now considering only points that lie on these $L$ lines, and call this set $P\subset{A\times{A}}$. Clearly, $|P|\approx{LN}$. Denote by $\Xi$ the set of slopes through the origin in this refined set $P$. More precisely, define $\Xi$ as follows:

$$\Xi:=\left\{\frac{b}{a}:(a,b)\in{P}\right\}.$$

The fact that $LN\leq{|A|^2}$ easily implies that

\begin{equation}
N\geq{\frac{M}{|A|^2}}.
\label{Nbound}
\end{equation}

Similarly, since $N\leq{A}$, it is clear that

\begin{equation}
L\geq{\frac{M}{|A|^2}}.
\label{Lbound}
\end{equation}

Let

$$A_x=\{y:(x,y)\in{P}\}$$

and

$$B_y=\{x:(x,y)\in{P}\}.$$

So, $A_x$ is the set of ordinates of $P$ for a fixed abscissa $x$, and $B_y$ is the set of abscissae of $P$ for a fixed ordinate $y$. 

A little more notation is required still. For some element $\xi\in{\Xi}$, let $P_{\xi}$ be the projection of points in $P$ on the line with equation $y=\xi{x}$ onto the x-axis. So,

$$P_{\xi}=\{x:(x,\xi{x})\in{P}\},$$

and since all lines with slope in $\Xi$ intersect $P$ with cardinality approximately $N$, it follows that $|P_{\xi}|\approx{N}$ for all $\xi{\in{\Xi}}.$ Note also, since $P\subset{A\times{A}}$, that $P_{\xi}$ is a subset of $A$.

The following lemma is taken from Rudnev(\cite{mishaSP}).

\begin{lemma}\label{theorem:misha} There exists a popular abscissa $x_0$ and a popular ordinate $y_0$ (popular here means that $|A_{x_0}|,|B_{y_0}|\gg{\frac{LN}{A}}$), as well as a large subset $\tilde{A}_{x_0}\subseteq{A_{x_0}}$, with

\begin{equation}
|\tilde{A}_{x_0}|\gg{\frac{LM}{|A|^3}},
\label{mishalemma1}
\end{equation}

such that for every $z\in{\tilde{A}_{x_0}}$,

\begin{equation}
|\tilde{A}_{\tilde{z}}:=P_{z/x_{0}}\cap{B_{y_0}}|\gg{\frac{LMN}{|A|^4}}.
\label{mishalemma2}
\end{equation}

\end{lemma}

Since the sum-product estimate is invariant under dilation, we may assume without loss of generality that $x_0=1$. This means that elements of $A_{x_0}$ are also the popular slopes described above, i.e. 

$$A_{x_0}\subset{\Xi}.$$

\subsection{Application of the covering lemma}

Since Lemma \ref{theorem:covering} is applied in a similar manner several times throughout the remainder of the proof, it is worthwhile highlighting this in advance, to avoid repetition.

\begin{application}\label{theorem:application} If $A'\subset{A}$ and $\xi\in{\Xi}$, then $\pm{\xi{A'}}$ can be 90\% covered by `at most' $\frac{K|A|}{N}$ translates of $A$.

\end{application}

\begin{proof} 90\% of $\xi{A'}$ can be covered by `at most'

$$\frac{|\xi{A'}+\xi{P_{\xi}}|}{|\xi{P_{\xi}}|}\ll{\frac{|A+A|}{N}}\leq{\frac{K|A|}{N}},$$

translates of $\xi{P_{\xi}}$, which is a subset of $A$. Similarly, 90\% of $-\xi{A'}$ can be covered by `at most'

$$\frac{|-\xi{A'}-\xi{P_{\xi}}|}{|\xi{P_{\xi}}|}\ll{\frac{|A+A|}{N}}\leq{\frac{K|A|}{N}},$$

translates of $\xi{P_{\xi}}\subset{A}$.

\end{proof}

\subsection{Four Cases}

For any set $B\subset{F}$, define $R(B)$ to be the set

$$R(B):=\left\{\frac{b_1-b_2}{b_3-b_4}:b_1,b_2,b_3,b_4\in{B},b_3\neq{b_4}\right\}.$$

The remainder of the proof is now divided into four cases corresponding to the nature of the sets $R(\tilde{A}_{x_0})$ and $R(B_{y_0})$.

\textbf{Case 1} $R(\tilde{A}_{x_0})$ is a subfield:

Case 1.1: First of all, it is possible that $|\tilde{A}_{x_0}|\leq{|A|^{6/7}}$. If so, then combining this bound with  \eqref{mishalemma1}, it can be deduced that

$$|A|^{27/7}\gg{LM}.$$

Using \eqref{Lbound}, rearranging, and then subsequently applying \eqref{dyadic} implies that

$$|A|^{41/7}\gg{M^2}\gg{\left(\frac{|A|^3}{K\log{|A|}}\right)^2},$$

and a simple rearrangement of the above gives

$$K\gg{\frac{|A|^{1/14}}{\log{|A|}}},$$

as required.

Case 1.2: The other possibility here is that $|\tilde{A}_{x_0}|>|A|^{6/7}$. In this case, since $R(\tilde{A}_{x_0})$ is a subfield,  Lemma \ref{theorem:KatzShen} tells us that $\tilde{A}_{x_0}\subseteq{cR(\tilde{A}_{x_0})+d}$ for some $c,d\in{F}$, and clearly

$$|A\cap{(cR(\tilde{A}_{x_0})+d)}|\geq{|\tilde{A}_{x_0}|}>|A|^{6/7}.$$
 
Then, the hypotheses of Theorem \ref{theorem:main} imply that

$$|R(\tilde{A}_{x_0})|^{14/25}\geq{|A\cap{(cR(\tilde{A}_{x_0})+d)}|}\geq{|\tilde{A}_{x_0}|},$$

and therefore

$$|R(\tilde{A}_{x_0})|\geq{|\tilde{A}_{x_0}|^{25/14}}.$$

For some $r\in{R(\tilde{A}_{x_0})}$, define $E_r(\tilde{A}_{x_0})$ to be the number of solutions to

$$a_1+ra_2=a_3+ra_4,$$

such that $a_1,a_2,a_3,a_4\in{\tilde{A}_{x_0}}$. A solution is considered to be trivial if $a_2=a_4$, and non-trivial otherwise. Summing over all $r\in{R(\tilde{A}_{x_0})}$ yields

\begin{align*}
\sum_{r\in{R(\tilde{A}_{x_0})}}E_r(\tilde{A}_{x_0})&=|\{\text{all trivial solutions}\}| + |\{\text{all non-trivial solutions}\}|
\\&\leq{|\tilde{A}_{x_0}|^{2}|R(\tilde{A}_{x_0})|+|\tilde{A}_{x_0}|^4}
\\&\leq{|\tilde{A}_{x_0}|^2|R(\tilde{A}_{x_0})|+|\tilde{A}_{x_0}|^{31/14}|R(\tilde{A}_{x_0})|}.
\end{align*}

Since the last term is dominant, it follows that

$$\sum_{r\in{R(\tilde{A}_{x_0})}}E_r(\tilde{A}_{x_0})\ll{|\tilde{A}_{x_0}|^{31/14}|R(\tilde{A}_{x_0})|}.$$

Hence, there exists some $r=\frac{p-q}{s-t}\in{R(\tilde{A}_{x_0})}$, such that $E_r(\tilde{A}_{x_0})\ll{|\tilde{A}_{x_0}|^{31/14}}$. Fix this $r$ and corresponding elements $p,q,s,t\in{\tilde{A}_{x_0}}$.

Let $\tilde{A}_{x_0}'$ be any subset of $\tilde{A}_{x_0}$ such that $|\tilde{A}_{x_0}'|\approx{|\tilde{A}_{x_0}|}$. This vague subset is introduced at this stage so that the covering lemma can later be applied effectively. Now, apply Cauchy-Schwarz in the usual way:

$$|\tilde{A}_{x_0}|^4\approx{|\tilde{A}_{x_0}'|^4}=\left(\sum_{x\in{\tilde{A}_{x_0}'+r\tilde{A}_{x_0}'}}\nu(x)\right)^2\leq{|\tilde{A}_{x_0}'+r\tilde{A}_{x_0}'||E_r(\tilde{A}_{x_0}')|}\ll{|\tilde{A}_{x_0}'+r\tilde{A}_{x_0}'||\tilde{A}_{x_0}|^{31/14}},$$

where $\nu(x)$ is the number of representations of $x$ as an element of $\tilde{A}_{x_0}'+r\tilde{A}_{x_0}'$. This implies that

$$|\tilde{A}_{x_0}|^{25/14}\ll{|\tilde{A}_{x_0}'+r\tilde{A}_{x_0}'|}\leq{|p\tilde{A}_{x_0}'-q\tilde{A}_{x_0}'+s\tilde{A}_{x_0}'-t\tilde{A}_{x_0}'|}.$$

Now, by Application \ref{theorem:application}, each of $p\tilde{A}_{x_0}',-q\tilde{A}_{x_0}',s\tilde{A}_{x_0}'$ and $-t\tilde{A}_{x_0}'$ can be 90\% covered by $\ll{\frac{K|A|}{N}}$ translates of $A$. Therefore, the subset $\tilde{A}_{x_0}'$ is chosen earlier in the proof so that each of $p\tilde{A}_{x_0}',-q\tilde{A}_{x_0}',s\tilde{A}_{x_0}'$ and $-t\tilde{A}_{x_0}'$ are covered completely by these copies of $A$. Applying the covering lemma four times yields

$$|\tilde{A}_{x_0}|^{25/14}\ll{|A+A+A+A|\left(\frac{K|A|}{N}\right)^4}\ll{\frac{K^{7}|A|^5}{N^4}},$$

where the last step follows from \eqref{finaltouch}. Applying bound \eqref{mishalemma1} from Lemma \ref{theorem:misha} and rearranging gives

$$M^{50/14}N^{6/14}\ll{K^7|A|^5|A|^{75/14}}.$$

Next, use \eqref{Nbound} to obtain

$$M^4\ll{K^7|A|^5|A|^{75/14}|A|^{12/14}}=K^7|A|^{157/14}.$$

Finally, apply \eqref{dyadic} and rearrange to get

$$K\gg{\frac{|A|^{1/14}}{(\log{|A|})^{4/11}}},$$

which is slightly better than required.

\textbf{Case 2} - $R(\tilde{A}_{x_0})\neq{R(B_{y_0})}$:

Case 2.1: There is some element $r\in{R(\tilde{A}_{x_0})}$ such that $r\notin{R(B_{y_0})}$. Fix this $r=\frac{a_1-a_2}{a_3-a_4}$ and elements $a_1,a_2,a_3,a_4\in{\tilde{A}_{x_0}}$ representing it. Since $r\notin{R(B_{y_0})}$, there exist only trivial solutions to

\begin{equation}
b_1+rb_2=b_3+rb_4,
\label{cseq}
\end{equation}

such that $b_1,b_2,b_3,b_4\in{B_{y_0}}$. Let $B_{y_0}'$ be some subset of $B_{y_0}$, with cardinality $|B_{y_0}'|\approx{|B_{y_0}|}$. Once again, this subset is required for the benefit of applying the covering lemma, and can be specified later. The absence of non-trivial solutions to \eqref{cseq} implies that 

$$|B_{y_0}'|^2\leq{|B_{y_0}'+rB_{y_0}'|}.$$

So,

$$\left(\frac{NL}{|A|}\right)^2\leq{|a_1B_{y_0}'-a_2B_{y_0}'+a_3B_{y_0}'-a_4B_{y_0}'|}.$$

By Application \ref{theorem:application}, each of $a_1B_{y_0},-a_2B_{y_0},a_3B_{y_0}$ and $-a_4B_{y_0}$ can be 90\% covered by $\ll{\frac{K|A|}{N}}$ translates of $A$. By choosing an appropriate subset $B_{y_0}'$, we can ensure that each of the four terms in the above sumset get fully covered. Therefore, Application \ref{theorem:application} is used four times in order to deduce that

$$\left(\frac{NL}{|A|}\right)^2\ll{\left(\frac{K|A|}{N}\right)^4|A+A+A+A|}\ll{\left(\frac{K|A|}{N}\right)^4K^3|A|}.$$

The above inequality can be rearranged into the form

$$M^2N^2\ll{K^7|A|^7}.$$

An application of \eqref{Nbound} and then subsequently \eqref{dyadic} implies that

$$K\gg{\frac{|A|^{1/11}}{(\log{|A|})^{4/11}}},$$

which is a better bound than required.

Case 2.2: There is some $r\in{R(B_{y_0})}$ such that $r\notin{R(\tilde{A}_{x_0})}$. Fix this $r=\frac{p-q}{s-t}$ as well as elements $p,q,s,t\in{B_{y_0}}$ representing $r$. Since $(p,y_0)\in{P}$ and all lines through the origin in $P$ support approximately $N$ points in $P$, it can be deduced that the line with gradient $\frac{y_0}{p}$ supports $\approx{N}$ points in $P$. The same can be said of the lines with gradients $\frac{y_0}{q}, \frac{y_0}{s}$ and $\frac{y_0}{t}$ respectively. In other words, each of $\frac{y_0}{p}, \frac{y_0}{q}, \frac{y_0}{s}$ and $\frac{y_0}{t}$ belong to $\Xi$.

Next, since $P$ is symmetric through the line $y=x$, it can be observed that $\frac{p}{y_0},\frac{q}{y_0},\frac{s}{y_0}$ and $\frac{t}{y_0}$ are also elements of $\Xi$, the set of slopes supporting $\approx{N}$ points from $P$. $r$ can be rewritten as

$$r=\frac{\frac{p}{y_0}-\frac{q}{y_0}}{\frac{s}{y_0}-\frac{t}{y_0}}.$$

Let $\tilde{A}_{x_0}'$ be a positively proportioned subset of $\tilde{A}_{x_0}$, to be chosen later in order to apply the covering lemma. Now, since $r\notin{R(\tilde{A}_{x_0}')}$, there exist only trivial solutions to the equation

$$a_1+ra_2=a_3+ra_4,$$

such that $a_1,a_2,a_3,a_4\in{\tilde{A}_{x_0}'}$. Therefore,

\begin{align*}
|\tilde{A}_{x_0}|^2\approx{|\tilde{A}_{x_0}'|^2}&\leq{|\tilde{A}_{x_0}'+r\tilde{A}_{x_0}'|}
\\&\leq{\left|\frac{p}{y_0}\tilde{A}_{x_0}'-\frac{q}{y_0}\tilde{A}_{x_0}'+\frac{s}{y_0}\tilde{A}_{x_0}'-\frac{t}{y_0}\tilde{A}_{x_0}'\right|}.
\end{align*}

By Application \ref{theorem:application}, each of $\frac{p}{y_0}\tilde{A}_{x_0},-\frac{q}{y_0}\tilde{A}_{x_0},\frac{s}{y_0}\tilde{A}_{x_0}$ and $-\frac{t}{y_0}\tilde{A}_{x_0}$ can be 90\% covered by $\ll{\frac{K|A|}{N}}$ translates of $A$. The subset $\tilde{A}_{x_0}'$ may be chosen earlier so that each of $\frac{p}{y_0}\tilde{A}_{x_0}',-\frac{q}{y_0}\tilde{A}_{x_0}',\frac{s}{y_0}\tilde{A}_{x_0}'$ and $-\frac{t}{y_0}\tilde{A}_{x_0}'$ are covered in their entirety by the translates of $A$.

Therefore, four applications of the covering lemma, along with \eqref{mishalemma1} and \eqref{finaltouch} yield

$$\left(\frac{LM}{|A|^3}\right)^2\ll{\left(\frac{K|A|}{N}\right)^4|A+A+A+A|}\ll{\frac{K^7|A|^5}{N^4}}.$$

This can be rearranged to give

$$M^4\ll{K^7|A|^{11}}.$$

Finally, apply \eqref{dyadic} to deduce that

$$K\gg{\frac{|A|^{1/11}}{(\log{|A|})^{4/11}}}.$$

From this point forward, we may assume that $R(\tilde{A}_{x_0})=R(B_{y_0})$.

\textbf{Case 3(Worst case)} - $R(\tilde{A}_{x_0})R(\tilde{A}_{x_0})\nsubseteq{R(\tilde{A}_{x_0})}=R(B_{y_0})$:

So, there must exist some $a_1,a_2,b_1,b_2,c_1,c_2,d_1,d_2\in{\tilde{A}_{x_0}}$ such that

$$\frac{a_1-b_1}{c_1-d_1}\frac{a_2-b_2}{c_2-d_2}\notin{R(\tilde{A}_{x_0})}.$$

\textbf{Claim}: This implies that there exists some elements $a,b,c,d,e,f\in{\tilde{A}_{x_0}}\subseteq{A_{x_0}}$ such that

$$r:=\frac{(a-b)(c-d)}{e-f}\notin{R(\tilde{A}_{x_0})}=R(B_{y_0}).$$

Suppose not. Then certainly

$$\frac{p-q}{s-t}=a_1-b_1\frac{a_2-b_2}{c_2-d_2}\in{R(A_{x_0})},$$

for some elements $p,q,s,t\in{A_{x_0}}$. Therefore, 

$$\frac{1}{c_1-d_1}\frac{p-q}{s-t}\notin{R(A_{x_0})}.$$

Note also that $R(A_{x_0})$ is closed under reciprocation, and thus

$$c_1-d_1\frac{s-t}{p-q}\notin{R(A_{x_0})},$$

which proves the claim.

Consider subsets $\tilde{A}_{\tilde{a}},\tilde{A}_{\tilde{c}}$ and $\tilde{A}_{\tilde{e}}$ of $B_{y_0}$; recall that these subsets were defined in the statement of Lemma \ref{theorem:misha}. 

Furthermore, let $\tilde{A}_{\tilde{a}}'\subset{\tilde{A}_{\tilde{a}}},\tilde{A}_{\tilde{c}}'\subseteq{\tilde{A}_{\tilde{c}}}$ and $\tilde{A}_{\tilde{e}}'\subseteq{\tilde{A}_{\tilde{e}}}$ be subsets with cardinality $|\tilde{A}_{\tilde{a}}'|\approx{|\tilde{A}_{\tilde{a}}|},|\tilde{A}_{\tilde{c}}'|\approx{|\tilde{A}_{\tilde{c}}|}$, and $|\tilde{A}_{\tilde{e}}'|\approx{|\tilde{A}_{\tilde{e}}|}$. These subsets are to be specified later in order to apply the covering lemma.

Since $\tilde{A}_{\tilde{a}}'$ and $\tilde{A}_{\tilde{e}}'$ are subsets of $B_{y_0}$, there exist only trivial solutions to the equation

$$a_1+ra_2=a_3+ra_4,$$

such that $a_1,a_3\in{\tilde{A}_{\tilde{e}}'}$ and $a_2,a_4\in{\tilde{A}_{\tilde{a}}'}$, as otherwise $r$ would be an element of $R(B_{y_0})$. This implies that

\begin{equation}
|\tilde{A}_{\tilde{e}}'||\tilde{A}_{\tilde{a}}'|\leq{|\tilde{A}_{\tilde{e}}'+r\tilde{A}_{\tilde{a}}'|}.
\label{solutions3}
\end{equation}

On the other hand, Lemma \ref{theorem:SPlun} can be applied with $X=\frac{c-d}{e-f}\tilde{A}_{\tilde{c}}'$ to bound the right hand side of the above inequality as follows:

\begin{align*}|\tilde{A}_{\tilde{e}}'+r\tilde{A}_{\tilde{a}}'|&\leq{\frac{|\frac{c-d}{e-f}\tilde{A}_{\tilde{c}}'+\tilde{A}_{\tilde{e}}'||\tilde{A}_{\tilde{c}}'+(a-b)\tilde{A}_{\tilde{a}}'|}{|\tilde{A}_{\tilde{c}}'|}}
\\&\leq{\frac{|c\tilde{A}_{\tilde{c}}'-d\tilde{A}_{\tilde{c}}'+e\tilde{A}_{\tilde{e}}'-f\tilde{A}_{\tilde{e}}'||a\tilde{A}_{\tilde{a}}'-b\tilde{A}_{\tilde{a}}'+\tilde{A}_{\tilde{c}}'|}{|\tilde{A}_{\tilde{c}}'|}}.
\end{align*}

Observe also that crucially, $a\tilde{A}_{\tilde{a}}'\subset{aP_{a}}\subset{A}$, and similarly $c\tilde{A}_{\tilde{c}}'\subset{A}$ and $e\tilde{A}_{\tilde{e}}'$, meaning there is no need to apply the covering lemma for these terms of the sum. Also, by definition, $\tilde{A}_{\tilde{c}}'\subset{A}$.

After combining this knowledge with \eqref{solutions3}, rearranging and applying \eqref{mishalemma2} to bound the left hand side, it follows that

$$\left(\frac{LMN}{|A|^4}\right)^3\ll{|A-d\tilde{A}_{\tilde{c}}'+A-f\tilde{A}_{\tilde{e}}'||A-b\tilde{A}_{\tilde{a}}'+A|}.$$

Next we must apply the covering lemma. By Application \ref{theorem:application}, each of $-d\tilde{A}_{\tilde{c}},-f\tilde{A}_{\tilde{e}}$ and $-b\tilde{A}_{\tilde{a}}$ can be 90\% covered $\ll{\frac{K|A|}{N}}$ translates of $A$. Therefore, the earlier choices of $\tilde{A}_{\tilde{c}}',\tilde{A}_{\tilde{e}}'$ and $\tilde{A}_{\tilde{a}}'$ are made so that each of $-d\tilde{A}_{\tilde{c}}',-f\tilde{A}_{\tilde{e}}'$ and $-b\tilde{A}_{\tilde{a}}'$ are covered completely by these translates of $A$. 

Therefore, only three applications of the covering lemma are needed in order to deduce that

$$\left(\frac{LMN}{|A|^4}\right)^3\ll{|A+A+A+A||A+A+A|\left(\frac{K|A|}{N}\right)^3}.$$

Applying \eqref{finaltouch} and \eqref{finaltouch2}, we get

$$\left(\frac{LMN}{|A|^4}\right)^3\ll{K^3|A|K^2|A|\left(\frac{K|A|}{N}\right)^3},$$

which can be rearranged to give

$$M^6\ll{K^8|A|^{17}}.$$

Finally, an application of \eqref{dyadic} implies that

$$K\gg\frac{|A|^{1/14}}{(\log{|A|})^{6/14}}.$$

\textbf{Case 4} $R(\tilde{A}_{x_0})R(\tilde{A}_{x_0})=R(\tilde{A}_{x_0})$ and $R(\tilde{A}_{x_0})+R(\tilde{A}_{x_0})\nsubseteq{R(\tilde{A}_{x_0})}$:

So, for some $a_1,a_2,a_3,a_4,b_1,b_2,b_3,b_4\in{\tilde{A}_{x_0}}$,

$$\frac{a_1-a_2}{b_1-b_2}+\frac{a_3-a_4}{b_3-b_4}\notin{R(\tilde{A}_{x_0})}=R(B_{y_0}).$$

Combining this with the knowledge that $R(\tilde{A}_{x_0})R(\tilde{A}_{x_0})=R(\tilde{A}_{x_0})$, we can deduce that there exist elements $p,q,s,t\in{\tilde{A}_{x_0}}$ such that

$$r:=\frac{p-q}{s-t}+1=\frac{b_3-b_4}{a_3-a_4}\frac{a_1-a_2}{b_1-b_2}+1=\frac{b_3-b_4}{a_3-a_4}\left(\frac{a_1-a_2}{b_1-b_2}+\frac{a_3-a_4}{b_3-b_4}\right)\notin{R(\tilde{A}_{x_0})}.$$

Let $\tilde{A}_{\tilde{p}}$ be as defined in Lemma \ref{theorem:misha}. Identify two subsets of positive proportion, $\tilde{A}_{\tilde{p}}'\subset{\tilde{A}_{\tilde{p}}}$ and $B_{y_0}'\subset{B_{y_0}}$, to be specified later for the purpose of applying the covering lemma. 

Applying Lemma \ref{theorem:CPlun} with $X=B_{y_0}'$, there exists a further subset $B_{y_0}''\subseteq{B_{y_0}'}$, with $|B_{y_0}''|\approx{|B_{y_0}'|}$, such that

\begin{equation}
\left|B_{y_0}''+\tilde{A}_{\tilde{p}}'+\left(\frac{p-q}{s-t}\right)\tilde{A}_{\tilde{p}}'\right|\ll{\frac{|A+A|}{|B_{y_0}|}\left|B_{y_0}'+\frac{p-q}{s-t}\tilde{A}_{\tilde{p}'}\right|}
\label{plunnecke}
\end{equation}

Now, since $B_{y_0}''$ and $\tilde{A}_{\tilde{p}}'$ are subsets of $B_{y_0}$, there exist only trivial solutions to

$$x_1+rx_2=x_3+rx_4,$$

such that $x_1,x_3\in{B_{y_0}''}$ and $x_2,x_4\in{\tilde{A}_{\tilde{p}}'}$, as otherwise $r\in{R(B_{y_0})}$. This implies that

$$\left(\frac{LN}{|A|}\right)\left(\frac{LMN}{|A|^4}\right)\ll{|B_{y_0}''||\tilde{A}_{\tilde{p}}'|}\ll{|B_{y_0}''+r\tilde{A}_{\tilde{p}}'|},$$

where the leftmost inequality is a consequence of \eqref{mishalemma2} and the lower bound on $|B_{y_0}|$ established earlier in Lemma \ref{theorem:misha}.

Combining this inequality with \eqref{plunnecke} and rearranging gives

\begin{equation}
\left(\frac{LN}{|A|}\right)^2\left(\frac{LMN}{|A|^4}\right)\ll{K|A|\left|B_{y_0}'+\left(\frac{p-q}{s-t}\right)\tilde{A}_{\tilde{p}}'\right|}.
\label{messy}
\end{equation}

Clearly,

$$\left|B_{y_0}'+\left(\frac{p-q}{s-t}\right)\tilde{A}_{\tilde{p}}'\right|\ll{|sB_{y_0}'-tB_{y_0}'+p\tilde{A}_{\tilde{p}}'-q\tilde{A}_{\tilde{p}}'|}.$$

Note also that $p\tilde{A}_{\tilde{p}}\subseteq{pP_{p}}\subseteq{A}$. Therefore,

$$\left|B_{y_0}'+\left(\frac{p-q}{s-t}\right)\tilde{A}_{\tilde{p}}'\right|\ll{|sB_{y_0}'-tB_{y_0}'+A-q\tilde{A}_{\tilde{p}}'|}.$$

Finally, three applications of the covering lemma are required. By Application \ref{theorem:application}, $sB_{y_0}$ and $-tB_{y_0}$ can be 90\% covered by `at most' $\frac{K|A|}{N}$ translates of $A$. The earlier choice of $\tilde{A}_{\tilde{s}}'$ should be made so as to ensure that both $sB{y_0}'$ and $-tB_{y_0}'$ get fully covered by these translates of $A$. Similarly, $-q\tilde{A}_{\tilde{p}}'$ can be fully covered by $\ll{\frac{K|A|}{N}}$ translates of $A$. It follows that

$$\left|B_{y_0}'+\left(\frac{p-q}{s-t}\right)\tilde{A}_{\tilde{p}}'\right|\ll{|A+A+A+A|\left(\frac{K|A|}{N}\right)^3}\ll{\frac{K^6|A|^4}{N^3}}.$$

Combining the above with \eqref{messy} and rearranging gives

$$M^4\ll{K^7|A|^{11}}.$$

Finally, an application of \eqref{dyadic}, leads to the conclusion that

$$K\gg{\frac{|A|^{1/11}}{(\log{|A|})^{4/11}}}.$$
\begin{flushright}
\qedsymbol
\end{flushright}

\bibliographystyle{plain}
\bibliography{reviewbibliography}

\end{document}